\documentclass{elsarticle}

\usepackage{lineno,hyperref}
\modulolinenumbers[5]
\usepackage{ifdraft}
\usepackage{amssymb,amsmath,amscd,epsfig,amsthm}
\usepackage{dsfont,graphicx}
\usepackage{fullpage,verbatim} 
\usepackage{tikz}
\usetikzlibrary{arrows,automata,shapes,backgrounds,shapes.symbols}
\usetikzlibrary{positioning}
\usetikzlibrary{decorations.pathmorphing}
\usetikzlibrary{snakes}
\usetikzlibrary{trees,decorations}
\usetikzlibrary{calc}
\usetikzlibrary{shapes.callouts,shapes.arrows}
\tikzstyle{btn}=[minimum size=15pt,color=white,draw=black,fill=black,circle]
\tikzstyle{ntb}=[color=black,draw=black,fill=white,circle]

\journal{Journal of \LaTeX\ Templates}

\newcommand{\Aut}{\mathop{\rm Aut}\nolimits}
\newcommand{\St}{\mathop{\rm Stab}}

\newcommand{\G}{\mathcal G}
\newcommand{\Hh}{\mathcal H}
\newcommand{\A}{\mathcal A}
\newcommand{\B}{\mathcal B}
\newcommand{\C}{\mathcal C}
\newcommand{\D}{\mathcal D}
\newcommand{\Z}{\mathbb Z}
\newcommand{\Sym}{\mathop{\rm Sym}\nolimits}

\newcommand{\Orb}{\mathop{\rm Orb}\nolimits}

\newcommand{\mz}{\mathfrak m}
\newcommand{\dz}{\mathfrak d}

\newcommand{\one}{\mathds 1}
\newcommand{\zero}{\mathds O}

\newcommand{\cw}{\mathrm{t}}
\newcommand{\dbel}{\langle\hat{\mathcal B}_3\rangle}
\newcommand{\bel}{\langle\mathcal B_3\rangle}

\usepackage{xcolor}

\newtheorem{theorem}{Theorem}[section]
\newtheorem{corollary}[theorem]{Corollary}
\newtheorem{proposition}[theorem]{Proposition}
\newtheorem{lemma}[theorem]{Lemma}
\newtheorem{conjecture}{Conjecture}
\theoremstyle{definition}
\newtheorem{definition}{Definition}
\newtheorem{example}{Example}
\newtheorem{remark}[theorem]{Remark}

\begin{document}

\begin{frontmatter}

\title{Orbit automata as a new tool to attack\\ the order problem in automaton groups}

%
\author[liafa]{Ines Klimann\fnref{mealym}}
\ead{klimann@liafa.univ-paris-diderot.fr}

\author[liafa]{Matthieu Picantin\fnref{mealym}}
\fntext[mealym]{Partially supported by the french \emph{Agence Nationale pour la~Recherche},
through the Project MealyM ANR-JCJC-12-JS02-012-01}
\ead{picantin@liafa.univ-paris-diderot.fr}

\author[usf]{Dmytro Savchuk\fnref{mealym,support}}
\fntext[support]{Partially supported by the New Researcher Grant and the Proposal Enhancement Grant from USF Internal Awards Program}
\ead{savchuk@usf.edu}

\address[liafa]{Univ Paris Diderot, Sorbonne Paris Cit\'e, LIAFA,
UMR 7089 CNRS, B\^atiment Sophie Germain, F-75013 Paris, France}
\address[usf]{Department of Mathematics and Statistics,
University of South Florida, 4202 E Fowler Ave, Tampa, FL 33620-5700, USA}

\begin{abstract}
We introduce a new tool, called the orbit automaton, that describes the action of an automaton group~$G$
on the subtrees corresponding to the orbits of~$G$ on levels of the tree.
The connection between~$G$ and the groups generated by the orbit automata
is used to find elements of infinite order in certain automaton groups for which other methods failed to work.
\end{abstract}


\end{frontmatter}


\section*{Introduction}
Groups generated by automata were formally introduced in 1960's~\cite{glushkov:automata,horejs:automata}, but gained a significant attention after remarkable discoveries in 1970's and 1980's that the class of these groups contains counterexamples to several long-standing conjectures in group theory. The first such evidence came in 1972 with the construction by Aleshin of an infinite periodic group generated by two initial automata~\cite{aleshin} (the complete proof can be found in~\cite{merzlyakov:periodic,grigorch:burnside}). But the field truly started to thrive after works of Grigorchuk~\cite{grigorch:burnside,grigorch:degrees} that introduced new methods of self-similarity and length contraction, and provided simpler counterexamples to the general Burnside problem, and the first counterexamples to the Milnor's problem on growth in groups~\cite{milnor}. We will also mention work of  Gupta and Sidki~\cite{gupta_s:burnside} that brought to life another series of related examples of infinite finitely generated $p$-groups and introduced a very powerful language of rooted trees to the field.

The class of automaton groups is particularly interesting from the computational viewpoint. The internal structure and complexity of these groups make computations by hands quite complicated, and sometimes infeasible. Even though the word problem is decidable for the whole class, other general algorithmic problems including the conjugacy problem and the isomorphism problem are known to be undecidable in general~\cite{sunic_v:conjugacy_undecidable}. The order problem was recently shown to be undecidable in the classes of semigroups generated by automata~\cite{gillibert:finiteness14} and groups generated by asynchronous automata~\cite{belk_b:undecidability}. However, the beauty of this class lies in the plethora of partial methods solving many algorithmic problems in majority of cases. For example, it was shown recently that the conjugacy problem and the order problem are decidable in the group of all, so-called, bounded automata~\cite{bondarenko2_sz:conjugacy13}.

Two software packages (\verb"FR"~\cite{bartholdi:fr} and \verb"AutomGrp"~\cite{muntyan_s:automgrp}) for \verb"GAP"~\cite{GAP4} system have been developed to address the computational demand in automaton groups and semigroups. Many of partial methods implemented in these packages rely heavily on the contraction of the length of the words while one passes to the sections at the vertices of the tree on which the group acts. However, not all automaton groups possess this property. In particular, such contraction rarely happens in groups generated by the \textit{reversible automata}. While working with these groups available software often fails to produce definite answers. At the same time, additional structure of reversible automata allows us to prove certain general results about the groups in this class. For example, in~\cite{klimann:finiteness} it is proved that the finiteness problem is decidable in the class of groups generated by 2-state reversible automata, in~\cite{klimann_ps:3state} it is proved that infinite groups generated by connected 3-state reversible automata always contain elements of infinite order. Further, in~\cite{godin:torsion-free} it is shown that invertible reversible automata that have no bireversible component with any number of states generate infinite torsion-free semigroups. There are also several papers~\cite{aleshin:free,gl_mo:compl,vorobets:aleshin,vorobets:series_free,steinberg_vv:series_free,savchuk_v:free_prods} devoted to the realization of free groups and free products of groups as groups generated by automata. All of the automata in the constructed families are reversible and by now all known free non-abelian automaton groups are either generated by reversible automata, or build from such groups~\cite{dangeli_r:free_nonreversible}. The above works use the interplay between the automaton and its dual to obtain the description of the group. Note, that there are other realizations of free groups as \emph{subgroups} of automaton groups (i.e. when we drop the ``self-similarity'' condition requiring the sections of the elements of the group to remain in the group). Specifically, the first realization of a nonabelian free group as a group generated by automata was achieved by Brunner and Sidki~\cite{brunner_s:glnz}.

In this paper we introduce a new tool, called the \emph{orbit automaton}, that can help to find the orders of some elements of automaton groups for which other methods fail to work. The main idea can be described as follows. Suppose an automaton group~$G$ acts not spherically-transitively on a rooted tree~$X^*$ consisting of all finite words over some finite alphabet~$X$, i.e. there is the smallest level $l$ of the tree on which $G$ acts non-transitively. Each orbit~$\mathcal O$ of the action of~$G$ on~$X^l$ in this case induces an invariant under~$G$ subtree~$T_{ \mathcal O}$ of~$X^*$ consisting of all words over~$X$ that have subwords only from~$\mathcal O$. This tree is not regular any more: it contains the $l$-th level $X^l$ of~$X^*$ and each vertex of~$X^l$ is the root of a regular $|\mathcal O|/|X^l|$-ary subtree of~$T_{\mathcal O}$. The orbit automaton associated with~$\mathcal O$ describes the action of~$G$ on these smaller degree subtrees. We show that for each automaton group~$G$ there is only finitely many different orbit automata (even when we consider the iterations of the above construction), and provide connections between~$G$ and groups generated by these orbit automata. Moreover, for reversible automata such connections provide more details about the groups~$\hat G$ generated by the automaton dual to the automaton generating~$G$.

As our main application we consider two groups that are generated by 4-state bireversible automata and prove that both of these groups contain torsion-free subsemigroups on at least two generators. These two groups were of a particular interest in the classification of all groups generated by 4-state automata over 2-letter alphabet started in~\cite{caponi:thesis2014}. Recall, that the classification of all 3-state automata over 2-letter alphabet was developed in~\cite{bondarenko_gkmnss:full_clas32_short}, where the notion of minimially symmetric automata was introduced. Namely, two automata are called \emph{symmetric} if it is possible to obtain one from the other by a sequence of symmetry operations that include permuting the states, permuting the letters of the alphabet, and passing to the inverse of the automaton. Two automata are called \emph{minimally symmetric} if their minimizations are symmetric. Both symmetry and minimal symmetry are obviously equivalence relations on the set of all automata that refine the relation induced by isomorphism of generated groups.

Among groups generated by 7471 non-minimally symmetric 4-state 2-letter automata all but 6 groups were shown in~\cite{caponi:thesis2014} either to be finite, or to contain an element of infinite order. Two out of these six groups are exactly groups studied here. The techniques developed in this paper potentially could be applied to handle the remaining four groups as well, as they are also generated by reversible automata.

The structure of the paper is as follows. In Section~\ref{sec:prelim} we set up the notation for automaton groups and recall necessary notions of the dual automaton, (bi)reversible automata, and the orbit tree. We introduce orbit automata and discuss the connections between the whole group and the groups generated by corresponding orbit automata in Section~\ref{sec:orbit_autom}. In Section~\ref{sec:applications} we present the main applications of the techniques developed in Section~\ref{sec:orbit_autom}.

\textbf{Acknowledgement.} The authors would like to thank Ievgen Bondarenko, Rostislav Grigorchuk, and Said Sidki for reading the first drafts of the paper and bringing up useful suggestions that enhanced the paper.

\section{Preliminaries}
\label{sec:prelim}

Let $X$ be a finite set, called an \emph{alphabet}, and let $X^*$ denote the set of all finite words over~$X$ (that can be though of as the free monoid generated by $X$). This set can be naturally endowed with a structure of a
rooted $|X|$-ary tree (where for a set $Y$ we denote by $|Y|$ its cardinality) by declaring that $v$ is adjacent to~$vx$ for
any $v\in X^*$ and $x\in X$. The empty word $\varepsilon$ corresponds to the root
of the tree and $X^n$ corresponds to the $n$-th level of the tree.
We will be interested in the groups of automorphisms and semigroups
of endomorphisms of~$X^*$ (as a graph). Any such endomorphism can be defined via
the notion of an initial automaton.

\begin{definition}
A \emph{Mealy automaton} (or simply \emph{automaton}) is a tuple
$\A=(Q,X,\pi,\lambda)$, where $Q$ is a set (the set of states), $X$ is a
finite alphabet, $\pi\colon Q\times X\to Q$ is a transition function
and $\lambda\colon Q\times X\to X$ is an output function. If the set
of states $Q$ is finite the automaton~$\A$ is called \emph{finite}. If
for every state $q\in Q$ the output function~$\lambda(q,\cdot)$ induces
a permutation of~$X$, the automaton~$\A$ is called invertible.
Selecting a state~$q\in Q$ produces an \emph{initial automaton}~$\A_q$.
\end{definition}

Automata are often represented by the \emph{Moore diagrams}. The
Moore diagram of an automaton $\A=(Q,X,\pi,\lambda)$ is a directed
graph in which the vertices are the states from~$Q$ and the labeled edges
have form
$q\xrightarrow{x|\lambda(q,x)}\pi(q,x)$
for~$q\in Q$ and~$x\in X$. Examples of Moore diagrams are shown in
Figure~\ref{fig:automH}.

Any initial automaton induces an endomorphism of the tree~$X^*$ (here we specifically view $X^*$ as a tree and not as a free monoid). Given a word
$v=x_1x_2x_3\ldots x_n\in X^*$, it scans its first letter $x_1$ and
outputs $\lambda(x_1)$. The rest of the word is handled in a similar
fashion by the initial automaton $\A_{\pi(x_1)}$. Formally speaking,
the functions $\pi$ and $\lambda$ can be extended recursively to~$\pi\colon
Q\times X^*\to Q$ and $\lambda\colon  Q\times X^*\to X^*$ via
\[\begin{array}{l}
\pi(q,x_1x_2\ldots x_n)=\pi(\pi(q,x_1),x_2x_3\ldots x_n),\\
\lambda(q,x_1x_2\ldots x_n)=\lambda(q,x_1)\lambda(\pi(q,x_1),x_2x_3\ldots x_n).\\
\end{array}
\]

By construction any initial automaton acts on~$X^*$ as an
endomorphism. In the case of invertible automaton it acts as an
automorphism. We will denote the group of all automorphisms of~$X^*$ by $\Aut(X^*)$.
For each invertible Mealy automaton $\A$ one can construct the \emph{inverse automaton} $\A^{-1}$ defined by swapping the labels of all arrows in the Moore diagram of~$\A$. It is easy to see that the states of~$\A^{-1}$ define the inverse transformations to those defined by corresponding states of~$\A$.

\begin{definition}
Let $\A$ be an (invertible) automaton over an alphabet $\Sigma$. The semigroup $\langle\A\rangle_+$ (group $\langle\A\rangle$) generated by all states of $\A$ viewed as endomorphisms (automorphisms) of the rooted tree $\Sigma^*$ under the operation of composition is called an \emph{automaton semigroup} (\emph{automaton group}).
\end{definition}

Another popular name for automaton groups and semigroups is
self-similar groups and semigroups
(see~\cite{nekrash:self-similar}).

Conversely, any endomorphism of~$X^*$ can be encoded by the action
of an initial automaton. In order to show this we need a notion of a
\emph{section} of an endomorphism at a vertex of the tree. Let $g$ be
an endomorphism of the tree~$X^*$ and $x\in X$. Then for any $v\in
X^*$ we have
\[g(xv)=g(x)v'\]
for some $v'\in X^*$. Then the map $g|_x\colon X^*\to X^*$ given by
\[g|_x(v)=v'\]
defines an endomorphism of~$X^*$ that is called the \emph{section} of
$g$ at vertex $x$. Furthermore,  for any $x_1x_2\ldots x_n\in X^*$
we define \[g|_{x_1x_2\ldots x_n}=g|_{x_1}|_{x_2}\ldots|_{x_n}.\]

Given an endomorphism $g$ of~$X^*$ we construct an initial automaton
$\A(g)$ whose action on~$X^*$ coincides with that of~$g$ as follows.
The set of states of~$\A(g)$ is the set~$\{g|_v\colon  v\in X^*\}$
of different sections of~$g$ at the vertices of the tree. The
transition and output functions are defined by
\[\begin{array}{l}
\pi(g|_v,x)=g|_{vx},\\
\lambda(g|_v,x)=g|_v(x).
\end{array}\]

Throughout the paper we will use the following convention. If $g$
and $h$ are the elements of some (semi)group acting on set $A$ and
$a\in A$, then
\begin{equation}
\label{eqn_conv}
gh(a)=h(g(a)).
\end{equation}

Taking into account convention~\eqref{eqn_conv} one can compute
sections of any element of an automaton semigroup as follows.
For~$g=g_1g_2\cdots g_n$ and $v\in X^*$, we have

\begin{equation}
\label{eqn_sections} g|_v=g_1|_v\cdot g_2|_{g_1(v)}\cdots
g_n|_{g_1g_2\cdots g_{n-1}(v)}.
\end{equation}

For any automaton group~$G$ there is a natural embedding
\[G\hookrightarrow G \wr \Sym(X)\]
of~$G$ into the permutational wreath product of~$G$ and the symmetric group~$\Sym(X)$ on~$X$ defined by
\[G\ni g\mapsto (g_1,g_2,\ldots,g_{|X|})\lambda(g)\in G\wr \Sym(X),\]
where $g_1,g_2,\ldots,g_{|X|}$ are the sections of~$g$ at the vertices
of the first level of~$X^*$, and $\lambda(g)$ is a permutation of~$X$ induced by the action of~$g$ on the first level of the tree.

The above embedding is convenient in computations involving the
sections of automorphisms, as well as for defining automaton groups. Sometimes the list of images of generators of an automaton group under this embedding is called the \emph{wreath recursion} defining the group.

In Section~\ref{sec:applications} we will need to work with dual automata, so we recall here necessary definitions and results. For any finite Mealy automaton one can construct the \emph{dual automaton}
defined by switching the roles of the set of states and the alphabet as well as
switching the transition and the output functions.

\begin{definition}
Given a finite automaton $\A=(Q,X,\pi,\lambda)$ its \emph{dual
automaton} $\hat\A$ is the finite automaton
$(X,Q,\hat\lambda,\hat\pi)$, where
\[\begin{array}{l}
\hat\lambda(x,q)=\lambda(q,x),\\
\hat\pi(x,q)=\pi(q,x)
\end{array}\]
for any $x\in X$ and $q\in Q$.
\end{definition}

Note that the dual of the dual of an automaton $\A$ coincides with~$\A$.
The semigroup~$\langle \hat\A\rangle_+$ generated by dual automaton~$\hat\A$
of automaton~$\A$ acts on the free monoid~$Q^*$. This
action induces the action on~$\langle \A\rangle_+$. Similarly, $\langle\A\rangle_+$ acts on~$\langle\hat\A\rangle_+$.

For an automaton (semi)group~$G$ generated by automaton $\A$, with a slight abuse of notation (since there might be many automata that define $G$), we will denote by $\hat G$ the (semi)group generated by the dual automaton $\hat\A$. In the case when the automaton generating $G$ is clear from the context we will call $\hat G$ the \emph{dual to~$G$ (semi)group}.

A particularly important class of automata is the class of
reversible automata.

\begin{definition}
An automaton $\A$ is called \emph{reversible} if its dual is invertible. If an automaton is invertible, reversible, and its inverse is reversible, it is called \emph{bireversible}.
\end{definition}

In particular, for any group generated by an invertible reversible automaton~$\A$,
one can consider the dual group generated by the dual automaton $\hat\A$.

The following proposition is proved in~\cite{vorobets:aleshin} by induction on level of the tree. With a slight abuse of notations we will denote by the same symbol the element of a free monoid and its
image under canonical epimorphism onto the corresponding semigroup.

\begin{proposition}
\label{prop_dual_sections}
Let $G$ be an automaton semigroup acting on~$X^*$ and generated by
the finite set $S$.  And let $\hat G$ be a dual semigroup to~$G$
acting on~$S^*$. Then for any $g\in G$ and $v\in X^*$ we have
$g|_v=v(g)$ in~$G$. Similarly, for any $g\in S^*$ and $v\in \hat G$,
$v|_g=g(v)$ in~$\hat G$.
\end{proposition}

Further, we will need a proposition relating the (semi)group generated by an automaton and the (semi)group generated by its dual.

\begin{proposition}[\cite{nekrash:self-similar,savchuk_v:free_prods,akhavi_klmp:finiteness_problem}]
\label{prop:dual_finite}
The (semi)group generated by an automaton $\A$ is finite if and only if the (semi)group generated by the automaton $\hat\A$ dual to~$\A$ is finite.
\end{proposition}

The last basic tool that we will need is the orbit tree of the action of an automaton group~$G$. Orbit trees in various forms have been studied earlier (see, for example,~\cite{gawron_ns:conjugation,bondarenko_s:sushch,klimann:finiteness,klimann_ps:3state,godin:torsion-free}) and describe the partition of the action of a group acting on a rooted tree into transitive components.

\begin{definition}
Let $G$ be an automaton group acting on a regular $|X|$-ary tree~$X^{*}$. The
\emph{orbit tree} of~$G$ is a graph whose vertices are the orbits of
$G$ on the levels of~$X^{*}$ and two orbits are adjacent if and only
if they contain vertices that are adjacent in~$X^{*}$. If, additionally, we label each edge connecting orbit $\mathcal O_1$ to orbit $\mathcal O_2$ on the next level by~$|\mathcal O_2|/|\mathcal O_1|$, we obtain a labeled graph that we call a \emph{labeled orbit tree}.
\end{definition}

Note that for each vertex in a labeled orbit tree of~$G$ the sum of the labels of all edges going down from this vertex always equals to the degree~$|X|$ of the tree~$X^*$. We also remark that one can define orbit tree of~$G$ in terms of connected components of powers of the dual of the automaton generating~$G$ (see~\cite{klimann_ps:3state}).

\section{Orbit Automaton}
\label{sec:orbit_autom}

We describe below a general construction that sometimes gives additional information about the structure of a self-similar group that acts non spherically transitively on the levels of the tree. First, we will need the following auxiliary definition.

\begin{definition}
Let $G$ be a automaton group generated by an automaton
over an alphabet~$X$.
The \emph{maximum transitivity level} $\cw(G)\in\mathbb N\cup\{0,\infty\}$ of the action of~$G$ on~$X^*$
is the maximum level of~$X^*$ on which $G$ acts transitively.
\end{definition}

Note that for an automaton $\A$ the maximum transitivity level of~$\langle\A\rangle$ coincides with the connection degree~$\curlywedge(\hat\A)$ of the dual automaton~$\hat\A$ defined in~\cite{klimann:finiteness}, see also~\cite{klimann_ps:3state}. The maximum transitivity level of~$G$ also has a natural interpretation as the length of the initial segment in the orbit tree before the first split.

Suppose that an automaton group~$G$ acts on~$X^*$ not spherically transitively, hence~$\cw(G)<\infty$. Then it follows from~\cite[Proposition 15]{klimann_ps:3state} (whose proof works for arbitrary number of states in the automaton) and Proposition~\ref{prop:dual_finite} that if all orbits of~$G$ on the level~$\cw(G)+1$ have size~$|X|^{\cw(G)}$ (\emph{i.e.} there are $|X|$ edges in the orbit tree going down from the vertex corresponding to the unique orbit on the level $X^{\cw(G)}$ and all of them are labeled by 1), then $G$ must be finite. Therefore,  the most interesting case
is when there is at least one orbit of size greater than $|X|^{\cw(G)}$.

\begin{definition}
\label{def:orbital_tree}
Let $G$ be an automaton group generated by an automaton over an alphabet~$X$ with $\cw(G)<\infty$ and let $\mathcal O$ be an orbit of~$G$ on the level~$\cw(G)+1$ of~$X^*$. Define the \emph{orbital tree~$T_{\mathcal O}$ of~$\mathcal O$} to be the subtree of~$X^*$
consisting of words whose all length~\((\cw(G)+1)\) subwords belong to~\(\mathcal{O}\).
\end{definition}

\begin{lemma}
\label{lem:structureTO}
Under the condition of Definition~\ref{def:orbital_tree} of the orbital tree~$T_{\mathcal O}$
the following holds.
\begin{itemize}
\item[(a)] $T_{\mathcal O}$ is an infinite spherically homogeneous tree whose vertices of levels 0 through $\cw(G)-1$ have $|X|$ children and all other vertices have $|\mathcal O|/|X|^{\cw(G)}$ children.
\item[(b)] $T_{\mathcal O}$ is invariant under the action of~$G$.
\end{itemize}
\end{lemma}

\begin{proof}
By definition no word of length less than $\cw(G)+1$ can contain a subword from~$X^{\cw(G)+1}\smallsetminus\mathcal O$, therefore all vertices of~$X^*$ up to level $\cw(G)$ are in~$T_{\mathcal O}$. Now suppose that $v=x_1x_2\ldots x_n$ with $n\geq\cw(G)$ is a vertex of~$T_{\mathcal O}$. Then $v'=x_{n-\cw(G)+1}x_{n-\cw(G)+2}\ldots x_n$ is a word from~$X^{\cw(G)}$, so there are exactly $k:=|\mathcal O|/|X|^{\cw(G)}$ letters $y_1,y_2,\ldots,y_k$ such that $v'y_i\in\mathcal O$ for all $1\leq i\leq k$. Therefore, the words $vy_i$, $1\leq i\leq k$ represent exactly $k$ children of vertex $v$ in~$T_{\mathcal O}$. This concludes the proof of item (a). We remark that it was crucial in the above argument that $\mathcal O$ is an orbit on level $\cw(G)+1$.

To show (b) we observe that if $v=x_1x_2\ldots x_n$ is a vertex of~$T_{\mathcal O}$, then for each $g\in G$ we have $g(v)=y_1y_2\ldots y_n$, where
\[y_iy_{i+1}\ldots y_{i+\cw(G)}=g|_{x_1x_2\ldots x_{i-1}}(x_ix_{i+1}\ldots x_{i+\cw(G)})\in\mathcal O\]
since $x_ix_{i+1}\ldots x_{i+\cw(G)}\in\mathcal O$ and $g|_{x_1x_2\ldots x_{i-1}}\in G$.
\end{proof}

Next, we introduce the notion of an \emph{orbit automaton}. The action of~$G$ on~$T_{\mathcal O}$ might not be faithful in general, so for a normal subgroup~$\St_{G}(T_{\mathcal O})$ of~$G$ consisting of elements of~$g$ that stabilize $T_{\mathcal O}$ we consider a group~$\overline{G}=G/\St_{G}(T_{\mathcal O})$ that acts faithfully on~$T_{\mathcal O}$. Throughout the paper for~$g\in G$ we denote by $\bar g$ the image of~$g$ under the canonical projection $G\to\overline G$. By Lemma~\ref{lem:structureTO} the orbital tree~$T_{\mathcal O}$ is naturally isomorphic to the tree~$X^{\cw(G)}Y^*$, where $Y=\{1,2,\ldots,k\}$ is an alphabet consisting of~$k=|\mathcal O|/|X|^{\cw(G)}$ letters. This isomorphism induces the action of~$\overline{G}$ on~$X^{\cw(G)}Y^*$ by automorphisms. The orbit automaton will define the action of sections of elements of~$\overline{G}$ on~$k$-ary subtrees of the form $vY^*$ for words $v$ of length at least $\cw(G)$.

More specifically, for each $v\in X^{\cw(G)}$ consider the $k$-ary subtree~$T_v$ of~$T_{\mathcal O}$ hanging down from vertex $v$.
There is a natural isomorphism
\[\psi_v\colon T_v\to Y^*\]
defined recursively by levels as follows. The root of~$T_v$ is sent by $\psi_v$ to the empty word in~$Y^*$. Suppose for~$w\in T_v$ the image $\psi_v(w)$ is defined, and $w$ has children $wx_1,wx_2,\ldots,wx_k$ in~$T_v$ ordered increasingly, then we define
\[\psi_v(wx_i)=\psi_v(w)i.\]

For each $g\in G$ the isomorphisms $\psi_v, v\in X^{\cw(G)}$ induce the action of~$\bar g\in\overline{G}$ on~$X^{\cw(G)}Y^*$ as follows. For~$v\in X^{\cw(G)}$ and $u\in Y^*$ define
\[\bar g(vu)=\bar g(v)\psi_{\bar g(v)}\bigl(\overline{g|_v}(\psi_v^{-1}(u))\bigr),\]
where $\bar g(v)$ is defined already since $v\in X^{\cw(G)}$, and $\overline{g|_v}(\psi_v^{-1}(u))$ is defined because $\psi_v^{-1}(u)\in T_{\mathcal O}$.

Therefore, we can naturally define the sections of elements of~$\overline{G}$ at vertices of~$X^{\cw(G)}$ that act on~$k$-ary trees isomorphic to~$Y^*$:
\[\bar g|_v(u)=\psi_{\bar g(v)}\bigl(\overline{g|_v}(\psi_v^{-1}(u))\bigr).\]
Since this automorphism of~$Y^*$ depends on~$g$ and  $v$, we will denote it as $(\overline{g|_v})_v$, where $\overline{g|_v}$ is an element of~$\overline{G}$.

\begin{theorem}
\label{thm:orbit_group}
Let $G$ be a group generated by automaton with the state set~$Q$ that does not act spherically transitively on the corresponding tree~$X^*$, and let $\mathcal O$ be an orbit of~$G$ on the level~$X^{\cw(G)+1}$ of the tree.

The group~$G_{\mathcal O}=\langle(\overline{g|_v})_v\colon g\in G,v\in X^{\cw(G)}\rangle<\Aut(Y^*)$ is an automaton group generated by an automaton with the state set
\[Q_{\mathcal O}=\{(\overline{q|_v})_v\colon q\in Q,v\in X^{\cw(G)}\}\]
and the wreath recursion defined as follows. Suppose $v=v_1v_2\ldots v_{\cw(G)}\in X^{\cw(G)}$ has children $vx_1,vx_2,\ldots,vx_k$ in~$T_{\mathcal O}$, where $k=|\mathcal O|/|X|^{\cw(G)}$. Then
\begin{itemize}
\item if $\cw(G)>0$:
\begin{equation}
\label{eqn:orbit_aut}
(\overline{q|_v})_v=\left(\bigl(\overline{(q|_{v_1})\bigl|_{v_2\ldots v_{\cw(G)}x_1}}\bigr)_{v_2\ldots v_{\cw(G)}x_1},\ldots,\bigl(\overline{(q|_{v_1})\bigl|_{v_2\ldots v_{\cw(G)}x_k}}\bigr)_{v_2\ldots v_{\cw(G)}x_k}\right)\sigma_{(\overline{q|_v})_v},
\end{equation}
\item if $\cw(G)=0$, then $v\in X^{\cw(G)}$ is an empty word $\varepsilon$ and
\begin{equation}
\label{eqn:orbit_aut2}
(\overline{q})_{\varepsilon}=(\overline{q|_{\varepsilon}})_{\varepsilon}=\left(\bigl(\overline{q|_{x_1}}\bigr)_{\varepsilon},\ldots,\bigl(\overline{q|_{x_k}}\bigr)_{\varepsilon}\right)\sigma_{(\overline{q})_{\varepsilon}},
\end{equation}
\end{itemize}
where $\sigma_{(\overline{q|_v})_v}(i)=j$ whenever $\bar q(vi)=v'j$ holds for~$i,j\in Y$, $v'\in X^{\cw(G)}$.
\end{theorem}

\begin{proof}
By construction for~$v=v_1v_2\ldots v_{\cw(G)}\in X^{\cw(G)}$, appending $i\in Y$ to~$v$ in~$X^{\cw(G)}Y^*$ corresponds to appending $x_i$ to~$v$ in~$X^*$. Therefore, we have
\[(\overline{g|_{v_1v_2\ldots v_{\cw(G)}}})_{v_1v_2\ldots v_{\cw(G)}}\bigl|_i=\bigl(\overline{(g|_{v_1})|_{v_2\ldots v_{\cw(G)}x_i}}\bigr)_{v_2\ldots v_{\cw(G)}x_i}.\]
This proves self-similarity of~$G_{\mathcal O}$ and equalities~\eqref{eqn:orbit_aut} and~\ref{eqn:orbit_aut2}.

Moreover, since
\begin{equation}
\label{eqn:prod}
\bigl(\overline{(gh)|_v}\bigr)_v=(\overline{g|_v})_v\cdot\bigl(\overline{h|_{g(v)}}\bigr)_{g(v)}
\end{equation}
and
\[\bigl(\overline{g|_v}\bigr)_v^{-1}=\bigl(\overline{g^{-1}|_{g(v)}}\bigr)_{g(v)},\]
the elements of~$Q_{\mathcal O}$ generate the whole $G_{\mathcal O}$.
\end{proof}

\begin{definition}
The group~$G_{\mathcal O}$ is called the \emph{orbit group} of~$G$ associated with the orbit~$\mathcal O$. The automaton~$\A_{\mathcal O}$ defined by wreath recursion~\eqref{eqn:orbit_aut} (or~\eqref{eqn:orbit_aut2} for~$\cw(G)=0$) and generating~$G_{\mathcal O}$ is called the \emph{orbit automaton} of~$G$ associated with the orbit~$\mathcal O$.
\end{definition}

\begin{remark}
\label{rem:infinite}
It might not be the case that each element of~$G_{\mathcal O}$ is a section of some element of~$\overline{G}$. In particular, it does not follow immediately that if $G_{\mathcal O}$ is infinite, then $\overline G$ is infinite. However, we do not have counterexamples to the following conjecture.
\end{remark}

\begin{conjecture}
For each finite automaton group $G$ all of the orbit groups of $G$ are finite.
\end{conjecture}

Despite the previous remark, we still can relate $\overline{G}$ and $G_{\mathcal O}$ via stabilizers of vertices from~$X^{\cw(G)}$ in~$\overline{G}$. Namely, for each $v\in X^{\cw(G)}$ the map
\begin{equation}
\label{eqn:hom_tau}
\begin{array}{lrcl}
\tau_v\colon &\St_{\overline{G}}(v)&\to &G_{\mathcal O},\\
&\bar g&\mapsto & (\overline{g|_v})_v.
\end{array}
\end{equation}
is a homomorphism since if $\bar g, \bar h\in\St_{\overline{G}}(v)$, then
\[\tau_v(gh)=((\overline{gh})|_v)_v=(\overline{g|_v})_v\cdot(\overline{h|_{g(v)}})_{g(v)}=(\overline{g|_v})_v\cdot(\overline{h|_v})_v=\tau_v(g)\tau_v(h).\]
This fact implies a particularly useful simple lemma that will be used later.

\begin{lemma}
\label{lem:tau}
If for some $v\in X^{\cw(G)}$ the subgroup~$\tau_v(\St_{\overline{G}}(v))$ of~$G_{\mathcal O}$ is infinite, then $\overline{G}$ is infinite.
\end{lemma}

\begin{proof}
The statement obviously holds since $\St_{\overline{G}}(v)$ is a subgroup of~$\overline{G}$.
\end{proof}

For each automaton group we can iterate the process of passing to the orbit groups only finite number of times. Therefore, this process allows us to associate to each automaton group a finite number of automaton groups acting on the regular trees of smaller degrees, thus ``decomposing'' the action of~$G$ into a finite number of smaller pieces.

\begin{proposition}
\begin{itemize}
\item[(a)] Each group which does not admit any nontrivial orbit groups is either finite, or acts spherically transitively on corresponding tree.

\item[(b)] For each automaton group~$G$ there is only finite number of automaton groups that can be obtained from~$G$ by iterative passing to the orbit groups.
\end{itemize}
\end{proposition}

\begin{proof}
First of all, if an automaton group does not have any nontrivial orbit groups, by definition this implies that either $\cw(G)=\infty$, in which case the group acts spherically transitively on the tree, or all of the orbits on the level $\cw(G)+1$ have the same size as the unique orbit on the level $\cw(G)$. In the latter case by~\cite[Proposition 15]{klimann_ps:3state} the group~$G$ must be finite.

To prove (b) it is enough to mention that the degree of the tree is reduced by at least one every time we pass to the orbit group. Therefore, since each automaton group has only finitely many orbit groups by definition, by iterative passing to orbit groups one can obtain only finitely many different automaton groups. Moreover, each such sequence will terminate in either a finite or spherically transitive group.
\end{proof}

We conclude this section with an observation saying that the orbit automata are sensitive to passing to symmetric automata because the order of the letters in the alphabet is used in the construction.

\begin{remark}
\label{rem:symmetric}
As shown in Example~\ref{ex:non_sym} below, it is possible for two symmetric automata $\A_1$ and $\A_2$ to generate isomorphic groups that have non minimally symmetric orbit automata and potentially non-isomorphic orbit groups (more precisely, $\langle\A_1\rangle$ could have an orbit automaton that is not minimally symmetric to any orbit automaton of~$\langle\A_2\rangle$).
\end{remark}

\begin{example}\label{ex-belleterra}
As a simple example, we will construct below an orbit automaton of the group $\dbel$ generated by the automaton $\hat{\mathcal B}_3$ with the state set~$X=\{\zero,\one\}$ dual to the Bellaterra automaton~$\mathcal B_3$ with the state set~$Q=\{a,b,c\}$ (both of these automata are depicted in Figure~\ref{fig:bellaterra}). Recall that both automata are bireversible and the group~$\bel$ generated by Bellaterra automaton~$\mathcal B_3$ is isomorphic to the free product~$C_2*C_2*C_2$ of three groups of order~$2$~\cite[p.25]{nekrash:self-similar}.

\begin{figure}[!h]
         \begin{center}
         \includegraphics{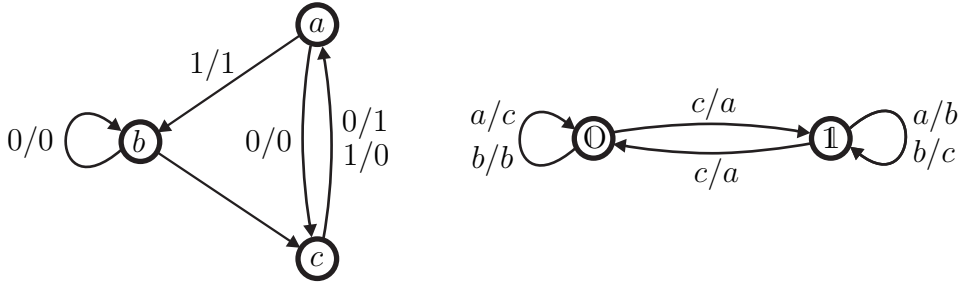}
         \end{center}
\caption{Automaton $\mathcal B_3$ generating the Bellaterra group (on left) and its dual $\hat{\mathcal B}_3$ (on right).\label{fig:bellaterra}}
\end{figure}

The wreath recursion of $\bel$ is:
\begin{equation}
\label{eqn_automH_def}
\begin{array}{lcl}
a&=&(c,b),\\
b&=&(b,c),\\
c&=&(a,a)\sigma.
\end{array}
\end{equation}
And the wreath recursion of $\dbel$ is:
\begin{equation}
\label{eqn_autom_dual_def}
\begin{array}{lcl}
\zero&=&(\zero,\zero,\one)(a\,c),\\
\one&=&(\one,\one,\zero)(a\,b\,c).\\
\end{array}
\end{equation}

The group $\dbel$ obviously acts transitively on the first level of $Q^*$, and on the second level it has two orbits: $\mathcal O_1=\{a^2,b^2,c^2\}$ and $\mathcal O_2=\{ab,ac,ba,bc,ca,cb\}$. Since $|\mathcal O_1|=|Q|$, the orbit automaton associated to~$\mathcal O_1$ is trivial. We will construct the orbit automaton associated to~$\mathcal O_2$.

The orbital tree $T_{\mathcal O_2}$ is shown in Figure~\ref{fig:orbit_tree} as a subtree of $Q^*$ (the edges of $T_{\mathcal O_2}$ are shown red bold). This is the same tree that was considered in~\cite{nekrash:self-similar} and consists of words over $Q$ which do not have equal consecutive letters.

\begin{figure}[!h]
        \begin{center}
        \includegraphics[width=420pt]{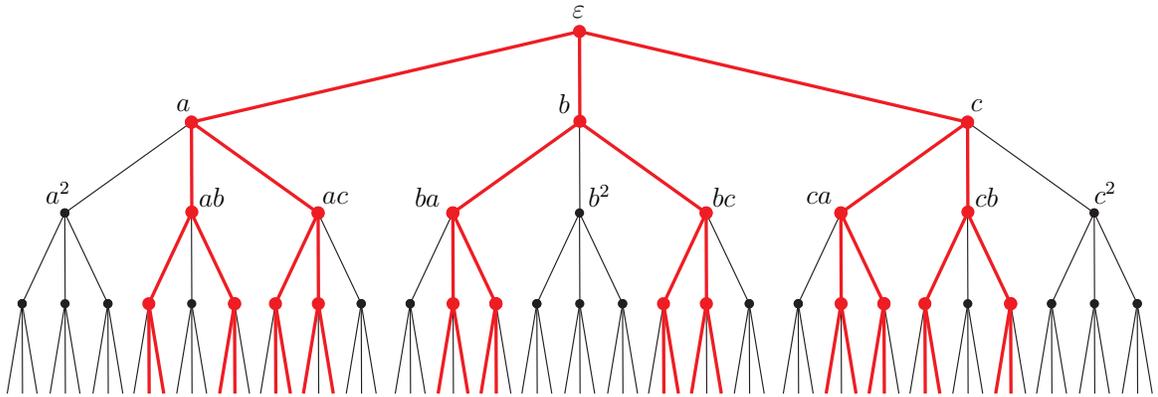}
        \end{center}
\caption{The orbital tree $T_{\mathcal O_2}$ (bold red) of the group generated by the dual to Bellaterra automaton $\hat{\mathcal B}_3$ associated to the orbit $\mathcal O_2=\{ab,ac,ba,bc,ca,cb\}$, shown as a subtree of the tree $Q^*$.\label{fig:orbit_tree}}
\end{figure}

By Theorem~\ref{thm:orbit_group} the orbit group $\dbel_{\mathcal O_2}$ acts on a binary tree $Y^*$ (for $Y=\{1,2\}$) and is generated by the set
\[Q_{\mathcal O_{2}}=\{(\overline{\zero|_a})_a,(\overline{\zero|_b})_b,(\overline{\zero|_c})_c,(\overline{\one|_a})_a,(\overline{\one|_b})_b,(\overline{\one|_c})_c\}
=\{\overline\zero_a,\overline\zero_b,\overline\one_c,\overline\one_a,\overline\one_b,\overline\zero_c\}.
\]
According to equation~\eqref{eqn:orbit_aut} the wreath recursion defining $\dbel_{\mathcal O_2}$ is as follows:
\begin{equation}
\label{eqn_autom_sectionsB}
\begin{array}{lclclcl}
\overline\zero_{a}&=&(\overline\zero_b,\overline\one_c)\sigma,&\qquad\qquad &\overline\one_{a}&=&(\overline\one_b,\overline\zero_c)\sigma,\\
\overline\zero_{b}&=&(\overline\zero_a,\overline\one_c)\sigma,&                  &\overline\one_{b}&=&(\overline\one_a,\overline\zero_c)\sigma, \\
\overline\zero_{c}&=&(\overline\zero_a,\overline\zero_b)\sigma,&             &\overline\one_{c}&=&(\overline\one_a,\overline\one_b).
\end{array}
\end{equation}
It is clear from above relations that $\zero_a=\zero_b$ and $\one_a=\one_b$. Therefore, the minimization of the orbit automaton $\mathcal A_{\mathcal O_2}$ is a 4-state automaton shown in Figure~\ref{fig:orb_autom_bel} that also generates $\dbel_{\mathcal O_2}$. Its wreath recursion is given below:
\begin{equation}
\label{eqn_autom_sectionsB_min}
\begin{array}{lcl}
\overline\zero_{a}&=&(\overline\zero_a,\overline\one_c)\sigma,\\ \overline\one_{a}&=&(\overline\one_a,\overline\zero_c)\sigma,\\
\overline\zero_{c}&=&(\overline\zero_a,\overline\zero_a)\sigma,\\           \overline\one_{c}&=&(\overline\one_a,\overline\one_a).
\end{array}
\end{equation}

\begin{figure}[!h]
        \begin{center}
        \includegraphics{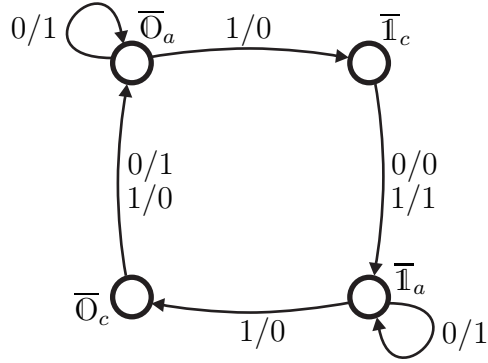}
        \end{center}
\caption{The minimization of the orbit automaton generating the orbit group $\dbel_{\mathcal O_2}$\label{fig:orb_autom_bel}.}
\end{figure}

The main part of the proof in~\cite{nekrash:self-similar} that $\bel$ has a structure of a free product is in showing that $\dbel$ acts spherically transitively on $T_{\mathcal O_2}$ proved in~\cite[Lemma~1.10.8]{nekrash:self-similar}. We suggest an alternative proof of this fact that uses the orbit group~$\dbel_{\mathcal O_2}$.

\begin{proposition}
The group $\dbel$ acts spherically transitively on $T_{\mathcal O_2}$.
\end{proposition}
\begin{proof}
First, we observe that the orbit group $\dbel_{\mathcal O_2}$ contains elements that generate cyclic groups acting spherically transitively on $Y^*$. For example, we can check using \verb"AutomGrp" package~\cite{muntyan_s:automgrp} that $\overline\zero_{a}\cdot\overline\one_{c}$ is one of such elements (for algorithms checking spherical transitivity we refer the reader to~\cite{bondarenko_gkmnss:full_clas32_short} and~\cite{steinberg:algorithmic}):
\begin{verbatim}
gap> D:=AutomatonGroup("0a=(0a,1c)(1,2),1a=(1a,0c)(1,2),0c=(0a,0a)(1,2),1c=(1a,1a)");
< 0a, 1a, 0c, 1c >
gap> IsSphericallyTransitive(0a*1c);
true
\end{verbatim}

Next, we observe that $\overline\zero_{a}\cdot\overline\one_{c}$ is in the image of the homomorphism $\tau_{a}$ defined in~\eqref{eqn:hom_tau}. Indeed, $\overline\zero^2$ fixes vertex~$a$ and by~\eqref{eqn:prod} we obtain
\[\tau_a(\overline\zero^2)=(\overline{\zero^2|_a})_a=(\overline{\zero}\cdot\overline{\one})_a=\overline{\zero}_a\cdot\overline{\one}_{\overline{\zero}(a)}=\overline{\zero}_a\cdot\overline{\one}_c.\]

Therefore, combining transitivity of $\dbel$ on the first level, and transitivity of $\overline\zero_{a}\cdot\overline\one_{c}$ on levels of $Y^*$, by induction on level we obtain transitivity of $\dbel$ on each level of $T_{\mathcal O_2}$.

\end{proof}

\end{example}

\section{Examples of applications}
\label{sec:applications}

We now turn our attention to two groups $\G$ and $\Hh$ generated by 4-state bireversible automata. We will prove that both of them contain torsion-free subsemigroups. The methods used for these groups are similar and use orbit trees of corresponding dual groups. After proving that $\Hh$ has a torsion-free subsemigroup, we use this fact to prove that $\G$ has a torsion-free subsemigroup by relating the group dual to~$\G$ to~$\Hh$ (thus we pass to the dual automaton twice: first from~$\G$ to~$\Hh$, and second to prove the result for~$\Hh$). As pointed out in the introduction, both of these groups are interesting, in particular, because automata that generate them represent 2 of 6 automata out of all 7421 non-symmetric 4-state invertible automata over 2-letter alphabet for which other ``standard'' methods of finding elements of infinite order did not work~\cite{caponi:thesis2014}. The relation between groups $\Hh$ and $\G$ is shown in Figure~\ref{fig:G_Hh}. First, from~$\Hh$ we pass to the dual group~$\hat\Hh$. The orbit group~$\hat\Hh_{\mathcal O_{qw}}$ of~$\hat\Hh$ is isomorphic to the group~$\G$. On the other hand, $\Hh$ is related to~$\G$ exactly in the same way as $\G$ is related to~$\Hh$.

\begin{figure}[!h]
\begin{center}
\[\begin{array}{cclcl}
\Hh&\leftrightsquigarrow&\hat\Hh&&\\
&&\hspace{0.6mm}\rotatebox[origin=c]{-90}{$\leadsto$}&&\\[0.5mm]
&&\hat\Hh_{\mathcal{O}_{qw}}\cong\G&\leftrightsquigarrow&\hat\G\\
&&&&\hspace{0.6mm}\rotatebox[origin=c]{-90}{$\leadsto$}\\[0.5mm]
&&&&\hat\G_{\mathcal{O}_{ac}}\cong\Hh
\end{array}\]
\caption{The relation between groups $\Hh$ and $\G$.\label{fig:G_Hh}}
\end{center}
\end{figure}

\subsection{Group~$\Hh$}
\label{ssec:autom_H}

Let $\Hh$ be a group defined by the following 4-state automaton $\B$ with the state set~$Q_{\Hh}=\{q,w,e,r\}$ depicted in the left side of Figure~\ref{fig:automH}. The wreath recursion defining $\Hh$ is given below:
\begin{equation}
\label{eqn_automH_def}
\begin{array}{lcl}
q&=&(w,w),\\
w&=&(e,q)\sigma,\\
e&=&(r,r)\sigma,\\
r&=&(q,e)\sigma.
\end{array}
\end{equation}

 \begin{figure}[!h]
         \begin{center}
         \includegraphics{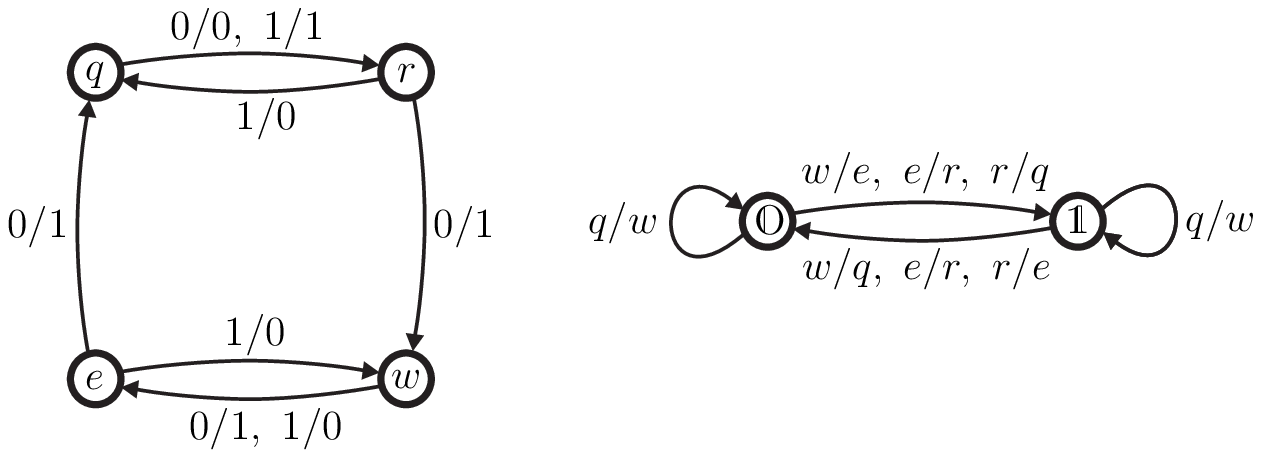}
         \end{center}
 \caption{Automaton $\mathcal B$ generating the group~$\Hh$ and its dual $\hat{\mathcal B}$.\label{fig:automH}}
 \end{figure}

\begin{proposition}
\label{prop:groupH}
Each finite word over $Q_{\Hh}$ that does not contain subwords
from~$\{q^2,w^2,e^2,r^2,qe,eq,wr,rw\}$ represents a nontrivial element in~$\Hh$.
In particular, the semigroup~$\langle er,qr\rangle_+$ is torsion-free.
\end{proposition}

Since the automaton $\B$ generating $\Hh$ is reversible, we can consider the dual group~$\hat\Hh$ generated by the automaton $\hat\B$ dual to~$\B$. The Moore diagram of this automaton is shown in the right side of Figure~\ref{fig:automH}, and its wreath recursion is as follows:

\begin{equation}
\label{eqn_autom_dual_def}
\begin{array}{lcl}
\zero&=&(\zero,\one,\one,\one)(q\,w\,e\,r),\\
\one&=&(\one,\zero,\zero,\zero)(q\,w)(e\,r).\\
\end{array}
\end{equation}

Let $T_{\hat\Hh}$ denote the labeled orbit tree of the action of~$\hat\Hh$ on the dual tree~$Q_{\Hh}^*$. The first six levels of are shown in Figure~\ref{fig:orbit_treeH}.
In particular, the maximum transitivity level~$\cw(\hat\Hh)$ equals $1$. Observe also that on the second level of $T_{\hat\Hh}$ we have two vertices that correspond to orbits of size~$8$.
According to~\cite[Lemmas 9 and 11]{klimann_ps:3state} this implies that there will be no edge labeled by~$3$ or~$4$ below the first level. We prove in Corollary~\ref{cor:path2} that there is a path in $T_{\hat\Hh}$ initiating at the root, whose edges below the first level are all labeled by~$2$ (as one can see from Figure~\ref{fig:orbit_treeH}, it is true up to level~$5$).

\begin{figure}[!ht]
\begin{center}
\input{otBTW.tex}
\end{center}
\caption{The labeled orbit tree~$T_{\hat\Hh}$ (up to level~$5$).}
\label{fig:orbit_treeH}
\end{figure}
Out of two orbits on the second level of of~$Q_{\Hh}^*$ (that are vertices on the second level of~$T_{\hat\Hh}$) we consider the orbit $\mathcal O_{qw}=\Orb_{\hat\Hh}(qw)$ consisting of the following vertices (words over $Q_{\Hh}$):
\[\mathcal O_{qw}=\{qw,wq,we,ew,er,re,qr,rq\}.\]

Below we will construct the orbit automaton of $\hat\Hh$ associated with $\mathcal O_{qw}$. Note, that the minimization of the orbit automaton of $\hat\Hh$ associated with the other orbit $\mathcal O_{q^2}$ of $\hat\Hh$ on the second level turns out to be trivial.

Consider the orbital tree~$T_{\mathcal O_{qw}}$ associated with the orbit~$\mathcal O_{qw}$. This subtree of~$Q_{\Hh}^*$ consists of finite words over~$Q_{\Hh}=\{q,w,e,r\}$ in which $q$ and~$e$ are always followed by~$w$ or~$r$, and $w$ and~$r$ are always followed by~$q$ or~$e$. By Lemma~\ref{lem:structureTO}(b) $T_{\mathcal O_{qw}}$ is invariant under the action of~$\hat\Hh$.

\begin{lemma}
\label{lem:stab_has_inf_indexH}
The stabilizer $K_{\hat\Hh}=\St_{\hat\Hh}(T_{\mathcal O_{qw}})$ of~$T_{\mathcal O_{qw}}$ in~$\hat\Hh$ has infinite index in~$\hat\Hh$.
\end{lemma}

\begin{proof}
By Lemma~\ref{lem:structureTO}(a) the root of the tree~$T_{\mathcal O_{qw}}$ has 4 children, and all other vertices have 2 children. Therefore the orbit group~$\hat\Hh_{\mathcal O_{qw}}$ associated with~$\mathcal O_{qw}$ acts on the binary tree~$Y^*$, where $Y=\{0,1\}$ is a 2-letter alphabet. By Theorem~\ref{thm:orbit_group} and according to the wreath recursion~\eqref{eqn_autom_dual_def}, $\hat\Hh_{\mathcal O_{qw}}$ is generated by
\begin{multline*}
Q_{\mathcal O_{qw}}=\{(\overline{\zero|_q})_q,(\overline{\zero|_w})_w,(\overline{\zero|_e})_e,(\overline{\zero|_r})_r,(\overline{\one|_q})_q,(\overline{\one|_w})_w,(\overline{\one|_e})_e,(\overline{\one|_r})_r\}
=\{\overline\zero_q,\overline\one_w,\overline\one_e,\overline\one_r,\overline\one_q,\overline\zero_w,\overline\zero_e,\overline\zero_r\}.
\end{multline*}

Further, according to~\eqref{eqn:orbit_aut} these generators satisfy the following wreath recursion:

\begin{equation}
\label{eqn_autom_sectionsH}
\begin{array}{lclclcl}
\overline\zero_{q}&=&(\overline\one_w,\overline\one_r)\sigma,&\qquad\qquad &\overline\one_{q}&=&(\overline\zero_w,\overline\zero_r),\\
\overline\zero_{w}&=&(\overline\zero_q,\overline\one_e),&                  &\overline\one_{w}&=&(\overline\one_q,\overline\zero_e), \\
\overline\zero_{e}&=&(\overline\one_w,\overline\one_r)\sigma,&             &\overline\one_{e}&=&(\overline\zero_w,\overline\zero_r),\\
\overline\zero_{r}&=&(\overline\zero_q,\overline\one_e),&                  &\overline\one_{r}&=&(\overline\one_q,\overline\zero_e).
\end{array}
\end{equation}

It is clear from the above relations that $\overline\zero_e=\overline\zero_q$, $\overline\zero_r=\overline\zero_w$, $\overline\one_e=\overline\one_q$, and $\overline\one_r=\overline\one_w$ hold. The minimization of the above automaton defines an automaton $\A$:

\begin{equation}
\label{eqn_autom_sectionsH_minimizedH}
\begin{array}{lcl}
\overline\zero_{q}&=&(\overline\one_w,\overline\one_w)\sigma,\\
\overline\zero_{w}&=&(\overline\zero_q,\overline\one_q),     \\
\overline\one_{q}&=&(\overline\zero_w,\overline\zero_w),\\
\overline\one_{w}&=&(\overline\one_q,\overline\zero_q).
\end{array}
\end{equation}

We would like to point out that the automaton $\A$ is obtained from the automaton~$\B$ by changing the permutations of all states. This is exactly the same kind of relation that was observed in Aleshin and Bellaterra families in~\cite{vorobets:series_free} and~\cite{steinberg_vv:series_free}.

By definition, the group~$\hat\Hh/K_{\hat\Hh}$ acts faithfully on the tree~$T_{\mathcal O_{qw}}$, which is isomorphic to the tree~$Q_{\Hh}\cdot Y^*$ (for~$Y=\{1,2\}$) via isomorphisms $\psi_s$, $s\in Q_{\Hh}$. This isomorphism between trees induces a faithful action of~$\hat\Hh/K_{\hat\Hh}$ on~$Q_{\Hh}\cdot Y^*$ defined by the following wreath recursion:
\begin{equation}
\label{eqn_action_subtreeH}
\begin{array}{lcl}
\overline{\zero}&=&(\overline\zero_q,\overline\one_w,\overline\one_q,\overline\one_w)(q\,w\,e\,r),\\
\overline{\one}&=&(\overline\one_q,\overline\zero_w,\overline\zero_q,\overline\zero_w)(q\,w)(e\,r),\\
\end{array}
\end{equation}
where the sections $\overline\zero_q$, $\overline\zero_w$, $\overline\one_q$, and $\overline\one_w$ are defined by~\eqref{eqn_autom_sectionsH_minimizedH}. We will identify $\hat\Hh/K_{\hat\Hh}$ with the group defined by this wreath recursion and will prove that it is infinite. This will imply that $K_{\hat\Hh}$ has infinite index in~$\hat\Hh$.

\medbreak The automaton $\mathcal A$ defined by~\eqref{eqn_autom_sectionsH_minimizedH} is bireversible, minimized, and its dual is a minimized 2-state automaton. Therefore, $\mathcal A$ is not $\mz\dz$-trivial and by~\cite{klimann:finiteness} it generates an infinite group~$\langle\mathcal A\rangle$. However, as pointed out in Remark~\ref{rem:infinite}, this does not imply immediately that $\hat\Hh/K_{\hat\Hh}$ is infinite as well: this might not be the case that for each $g\in\langle\overline\zero_q,\overline\zero_w,\overline\one_q,\overline\one_w\rangle=\langle\A\rangle$ there is an element of~$\hat\Hh/K_{\hat\Hh}$ whose section at some vertex of~$Q_{\Hh}\cdot Y^*$ is $g$.

Consider the stabilizer $S=\St_{\hat\Hh/K_{\hat\Hh}}(q)$ of vertex $q\in Q_{\Hh}\cdot Y^*$ in~$\hat\Hh/K_{\hat\Hh}$. It is a subgroup of finite index in~$\hat\Hh/K_{\hat\Hh}$ and by Reidemeister-Schreier procedure we find
\[S=\langle \overline{\one}\cdot\overline{\zero}^{-1}, \overline{\zero}\cdot\overline{\one}, \overline{\zero}^{-1}\cdot\overline{\one}\cdot\overline{\zero}^{-2}, \overline{\zero}^4, \overline{\zero}^2\cdot\overline{\one}\cdot\overline{\zero} \rangle.\]
As $S$ stabilizes vertex $q$, there is a natural homomorphism
\[\begin{array}{llll}
\tau_q\colon &S&\to &\langle\overline\zero_q,\overline\zero_w,\overline\one_q,\overline\one_w\rangle,\\
&h&\mapsto & h|_q.
\end{array}
\]
The image of~$S$ under $\tau_q$ is (taking into account that the generators of~$\langle\A\rangle$ are involutions):
\begin{multline*}
\label{eqn:defnS}
\tau_q(S)=\langle g_1=\overline\one_q\overline\zero_q,\quad g_2=\overline\zero_q\overline\zero_w,\quad g_3=\overline\one_w\overline\zero_w\overline\one_w\overline\zero_q, \quad
g_4=\overline\zero_q\overline\one_w\overline\one_q\overline\one_w,\quad g_5=(\overline\zero_q\overline\one_w)^2\rangle.
\end{multline*}
Observe that $\tau_q(S)$ is a normal subgroup of~$\langle\A\rangle=\langle\overline\zero_q,\overline\zero_w,\overline\one_q,\overline\one_w\rangle$. Indeed, since the generators of~$\langle\A\rangle$ all are involutions, it is enough to prove that conjugations by them normalizes $\tau_q(S)$. Below, using the fact that both $\langle\overline\zero_q,\overline\one_q\rangle$ and $\langle\overline\zero_w,\overline\one_w\rangle$ are isomorphic to the 4-element Klein group~$(\Z/2\Z)^2$, we show that this is indeed the case:
\[\begin{array}{llll}
g_1^{\overline\zero_q}=g_1,& g_1^{\overline\one_q}=g_1, &g_1^{\overline\zero_w}=g_5^{-1}g_4,& g_1^{\overline\one_w}=g_5^{-1}g_4,\\[1mm]
g_2^{\overline\zero_q}=g_3,& g_2^{\overline\one_q}=g_1g_3g_1, &g_2^{\overline\zero_w}=g_3,& g_2^{\overline\one_w}=g_5^{-1}g_2^{-1},\\[1mm]
g_3^{\overline\zero_q}=g_2,& g_3^{\overline\one_q}=g_1g_2g_1, &g_3^{\overline\zero_w}=g_2,& g_3^{\overline\one_w}=g_3g_5,\\[1mm]
g_4^{\overline\zero_q}=g_4^{-1},& g_4^{\overline\one_q}=(g_1g_4g_1)^{-1}, &g_4^{\overline\zero_w}=g_3g_4^{-1}g_2,& g_4^{\overline\one_w}=g_5^{-1}g_1,\\[1mm]
g_5^{\overline\zero_q}=g_5^{-1},& g_5^{\overline\one_q}=g_1g_5^{-1}g_1, &g_5^{\overline\zero_w}=g_3g_5^{-1}g_2,& g_5^{\overline\one_w}=g_5^{-1}.\\
\end{array}
\]

Moreover, since $g_2\overline\zero_w=\overline\zero_q$ and $g_1g_2\overline\zero_w=\overline\one_q$, we have
\[\langle\A\rangle=\langle\tau_q(S),\overline\zero_w,\overline\one_w\rangle.\]
Now, $\langle\overline\zero_w,\overline\one_w\rangle\cong (\Z/2\Z)^2$ implies
\[[\langle\A\rangle:\tau_q(S)]\leq 4.\]
But this means that $\tau_q(S)$ must be infinite as a subgroup of finite index in an infinite group~$\langle\A\rangle$ (recall that this is the case by non-$\mz\dz$-triviality according to~\cite{klimann:finiteness}). Finally, we conclude (as stated in Lemma~\ref{lem:tau}) that~$S$, and hence $\hat\Hh/K_{\hat\Hh}$ are also infinite.
\end{proof}

\begin{corollary}
\label{cor_stabs_diff}
The stabilizers of levels of~$T_{\mathcal O_{qw}}$ in~$\hat\Hh/K_{\hat\Hh}$ are pairwise different.
\end{corollary}
\begin{proof}
Since $\hat\Hh/K_{\hat\Hh}$ is infinite by Lemma~\ref{lem:stab_has_inf_indexH} and
all stabilizers of levels are finite index subgroups of~$\hat\Hh/K_{\hat\Hh}$,
they are all infinite. Let $g\in\hat\Hh$ be an arbitrary element such that $\bar g\in\hat\Hh/K_{\hat\Hh}$ is nontrivial and such that~$\bar g$ (and hence $g$)
fixes level~$n$ of~$T_{\mathcal O_{qw}}$. Let~$m\geq n+1$ be the smallest level on which $\overline{g}$
acts nontrivially. Then there exists a vertex $v=x_1x_2\ldots
x_{m-1}$ of~$T_{\mathcal O_{qw}}$, such that $\overline{g|_v}$ acts nontrivially on the
first level of corresponding binary tree. Then $\overline{g|_{x_1x_2\ldots x_{m-n-1}}}$ stabilizes the $n$-th level of~$T_{\mathcal O_{qw}}$ but does not stabilize the $(n+1)$-st. Note that by construction of~$T_{\mathcal O_{qw}}$ the word~$x_1x_2\ldots x_{m-n-1}$ also represents a vertex in~$T_{\mathcal O_{qw}}$.
\end{proof}

\begin{lemma}
\label{lem:transitivityH}
The action of~$\hat\Hh$ on the levels of~$T_{\mathcal O_{qw}}$ is transitive.
\end{lemma}

\begin{proof}
We prove it by induction on the level $n$. The statement is obviously true for~$n=1$. Assume that $\hat\Hh$ acts transitively on level $n$ for some $n\geq1$.

By Corollary~\ref{cor_stabs_diff} there is an
element $g\in\hat\Hh$ that fixes the $n$-th level of~$T_{\mathcal O_{qw}}$ and acts nontrivially on
the~$(n+1)$-st level. This means that there is a vertex $v$ on the $n$-th level of~$T_{\mathcal O_{qw}}$
such that $g(v)=v$ holds and $\overline{g|_v}$ acts nontrivially and, hence, transitively on the first level of corresponding binary tree.

Fix a vertex $vx$ of the $(n+1)$-st level of~$T_{\mathcal O_{qw}}$. By
induction assumption, for any vertex $y_1y_2\ldots y_{n+1}$ of the same level there is an element $h\in\hat\Hh$ that moves $v$ to
$y_1y_2\ldots y_n$. Then $g^kh$, where $k$ is 0 or 1, will move
$vx$ to~$y_1y_2\ldots y_{n+1}$. Thus, $\hat\Hh$ acts transitively on the levels of~$T_{\mathcal O_{qw}}$.
\end{proof}

The last lemma can be rephrased in terms of the orbit tree of $\hat\Hh$.

\begin{corollary}
\label{cor:path2}
The orbit tree of $\hat\Hh$ contains an infinite path from the root labeled by $42^{\infty}$.
\end{corollary}

\begin{proof}
By Lemma~\ref{lem:transitivityH} each level of~$T_{\mathcal O_{qw}}$ corresponds to one vertex in~$T_{\hat\Hh}$. Since $(i+1)$-st level of~$T_{\mathcal O}$ contains twice as many vertices as the $i$-th one for~$i\geq 1$, the edges in~$T_{\hat\Hh}$ connecting corresponding vertices are all labeled by~$2$. The path in the statement of this corollary is the image of~$T_{\mathcal O_{qw}}$ under the quotient map from~$Q_{\Hh}^*$ to~$T_{\hat\Hh}$.
\end{proof}

Finally, we provide the proof of Proposition~\ref{prop:groupH}.

\begin{proof}[Proof of Proposition~\ref{prop:groupH}.]
The elements of~$\Hh$ corresponding to vertices $(er)^ne$, $(eq)^nqr$, $n\geq0$ of odd and even levels of~$T_{\mathcal{O}_{qw}}$ respectively, act nontrivially on the first level. By Lemma~\ref{lem:transitivityH} each word $v$ corresponding to a vertex in~$T_{\mathcal O_{qw}}$ can be moved by $\hat\Hh$ to one of these vertices. But this means that $v$ has a nontrivial section and, thus, must be a nontrivial element of~$\Hh$ itself.
\end{proof}

\subsection{Group~$\G$}
\label{ssec:autom_G}

\begin{figure}[!h]
         \begin{center}
         \includegraphics{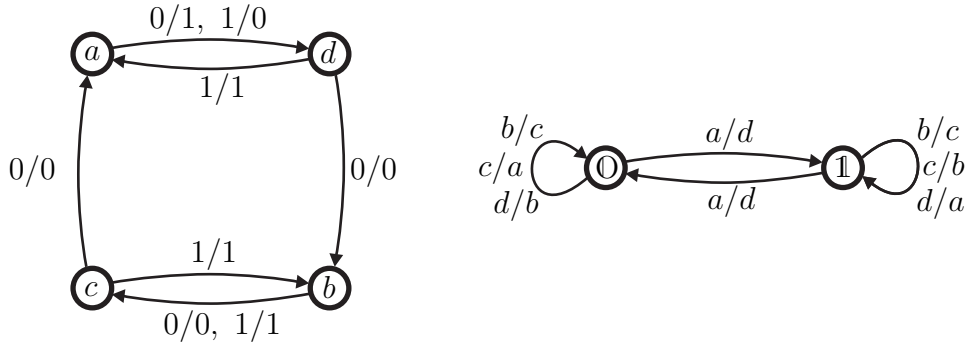}
         \end{center}
 \caption{Automaton $\mathcal C$ generating the group~$\G$ and its dual $\hat{\mathcal C}$.\label{fig:automG}}
 \end{figure}

We now turn our attention to the group~$\hat\Hh_{\mathcal O_{qw}}$ generated by the 4-state automaton constructed in the previous subsection. To simplify the notation here we will now denote this group by $\G$, and also we will rename the generators of~$\G$ (and the states of the generating automaton) according to the following rule: $a\leftrightarrow\overline\zero_q$, $b\leftrightarrow\overline\one_q$, $c\leftrightarrow\overline\zero_w$, $d\leftrightarrow\overline\one_w$. Finally, the automaton with the state set $Q_{\G}=\{a,b,c,d\}$ generating~$\G$ will be denoted by $\C$. This automaton is depicted in the left side of Figure~\ref{fig:automG}. Its defining wreath recursion is as follows:
\begin{equation}
\label{eqn_autom_def}
\begin{array}{lcl}
a&=&(d,d)\sigma,\\
b&=&(c,c),\\
c&=&(a,b),\\
d&=&(b,a).
\end{array}
\end{equation}

All generators of~$\G$ have order 2 and the subgroups $\langle a,b\rangle$ and $\langle c,d\rangle$ are both isomorphic to~$(\Z/2\Z)^2$.

\begin{figure}[!ht]
\begin{center}
\input{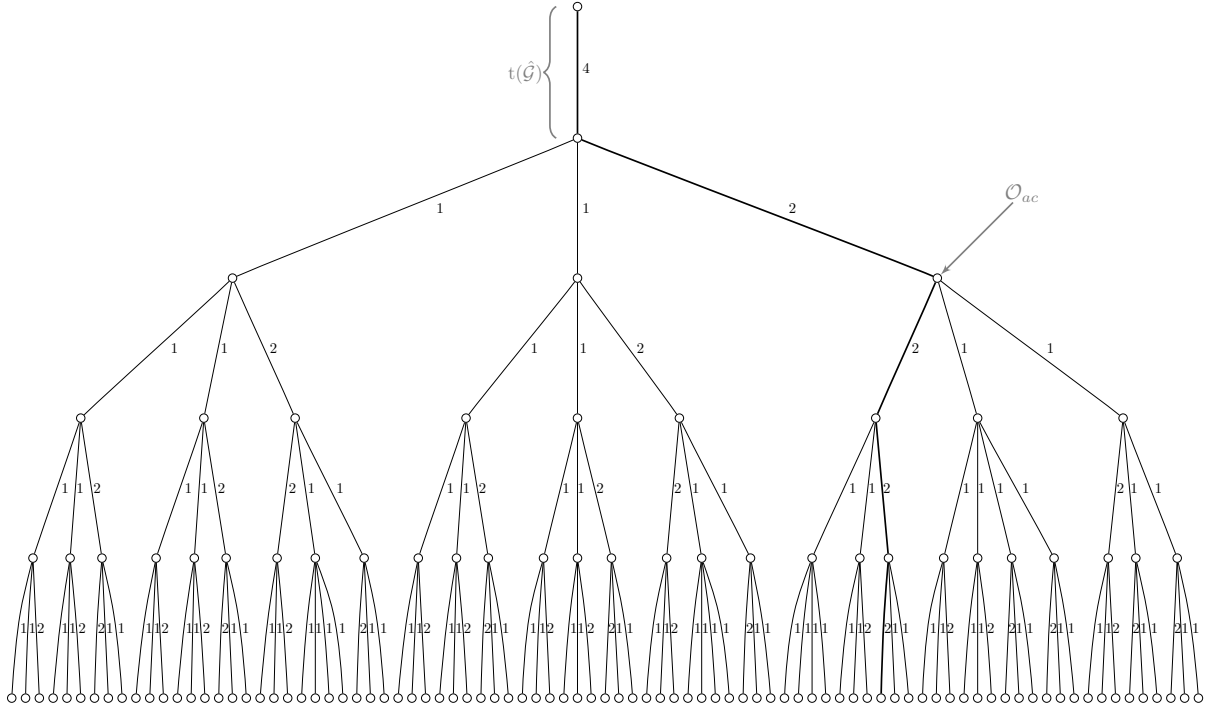}
\end{center}
\caption{The labeled orbit tree~$T_{\hat\G}$ (up to level~5).}
\label{fig:orbit_treeG}
\end{figure}

\begin{proposition}
\label{prop:groupG}
Each finite word over $Q_{\G}$ that does not contain subwords 
from~$\{a^2,b^2,c^2,d^2,ab,ba,cd,dc\}$ represents a nontrivial element in~$\G$.
In particular, the semigroups $\langle ac,ad,bc,bd\rangle_+$ and $\langle ca,da,cb,db\rangle_+$ are torsion-free.
\end{proposition}
The proof is very similar to the proof of Proposition~\ref{prop:groupH}, but by using this proposition one part of the proof becomes simpler. Since the automaton $\C$ generating $\G$ is reversible, we can consider the dual group~$\hat\G$ generated by the automaton $\hat\C$ dual to~$\C$, shown in the right side of Figure~\ref{fig:automG}. Its defining wreath recursion is given below:
\begin{equation}
\label{eqn_autom_dual_defG}
\begin{array}{lcl}
\zero&=&(\one,\zero,\zero,\zero)(a\,d\,b\,c),\\
\one&=&(\zero,\one,\one,\one)(a\,d)(b\,c).\\
\end{array}
\end{equation}
Let $T_{\hat\G}$ denote the labeled orbit tree of the action of~$\hat\G$ on the dual tree~$Q_{\G}^*$. The first six levels of are shown in Figure~\ref{fig:orbit_treeG}. Observe that on the second level we have two orbits of size~$4$ and one orbit of size~$8$. Again according to~\cite[Lemmas 9 and 11]{klimann_ps:3state} this implies that there will be no edge labeled by~$3$ or~$4$ below the first level of $T_{\hat\G}$ and that there will be at most one path from the root to infinity in~$T_{\hat\G}$ whose all edges below level one are labeled by~$2$. We will prove in Corollary~\ref{cor:path2G} that there is such a path (as one can see from Figure~\ref{fig:orbit_treeG}, it is true up to level~5).

In the case of~$\hat\G$ there is only one orbit $\mathcal O_{ac}=\Orb_{\hat\G}(ac)$ of size~$8$ on the second level of~$Q_{\G}^*$ (that is a vertex on the second level of~$T_{\hat\G}$). This orbit consists of the following vertices (words over $Q_{\G}$):
\[\mathcal O_{ac}=\{ac,ad,bc,bd,ca,da,cb,db\}.\]
By Lemma~\ref{lem:structureTO}(b), the orbital tree~$T_{\mathcal O_{ac}}$---defined as the subtree of~$Q_{\G}^*$ consisting of finite words over~$Q_{\G}=\{a,b,c,d\}$ in which $a$ and~$b$ are always followed by~$c$ or~$d$, and $c$ and~$d$ are always followed by~$a$ or~$b$---is invariant under the action of~$\hat\G$.

\begin{lemma}
\label{lem:stab_has_inf_indexG}
The stabilizer $K_{\hat\G}=\St_{\hat\G}(T_{\mathcal O_{ac}})$ of~$T_{\mathcal O_{ac}}$ in~$\hat\G$ has infinite index in~$\hat\G$.
\end{lemma}

\begin{proof}
Just as in the proof of Lemma~\ref{lem:stab_has_inf_indexH} by Lemma~\ref{lem:structureTO}(a) the root of the tree~$T_{\mathcal O_{ac}}$ has 4 children, and all other vertices have 2 children. Therefore, the orbit group~$\hat\G_{\mathcal O_{ac}}$ associated with~$\mathcal O_{ac}$ acts on the binary tree~$Y^*$, where $Y=\{0,1\}$ is a 2-letter alphabet. By Theorem~\ref{thm:orbit_group} and according to the wreath recursion~\eqref{eqn_autom_dual_defG},
$\hat\G_{\mathcal O_{ac}}$ is generated by
\begin{multline*}
Q_{\mathcal O_{ac}}=\{(\overline{\zero|_a})_a,(\overline{\zero|_b})_b,(\overline{\zero|_c})_c,(\overline{\zero|_d})_d,(\overline{\one|_a})_a,(\overline{\one|_b})_b,(\overline{\one|_c})_c,(\overline{\one|_d})_d\}
=\{\overline\one_a,\overline\zero_b,\overline\zero_c,\overline\zero_d,\overline\zero_a,\overline\one_b,\overline\one_c,\overline\one_d\}.
\end{multline*}

Further, according to~\eqref{eqn:orbit_aut} these generators satisfy the following wreath recursion (after taking into account that  $\overline\zero_a=\overline\zero_b$, $\overline\zero_c=\overline\zero_d$, $\overline\one_a=\overline\one_b$, and $\overline\one_c=\overline\one_d$ hold):

\begin{equation}
\label{eqn_autom_sectionsG}
\begin{array}{lcl}
\overline\zero_{a}&=&(\overline\zero_{c},\overline\zero_{c}),\\
\overline\zero_{c}&=&(\overline\one_{a},\overline\zero_{a})\sigma,\\
\overline\one_{a}&=&(\overline\one_{c},\overline\one_{c})\sigma,\\
\overline\one_{c}&=&(\overline\zero_{a},\overline\one_{a})\sigma.
\end{array}
\end{equation}

The key point of the proof is the observation that we get back exactly the automaton~$\B$ from the previous subsection defined by~\eqref{eqn_automH_def} under the state correspondence $q\leftrightarrow\overline\zero_a$, $w\leftrightarrow\overline\zero_c$, $e\leftrightarrow\overline\one_a$, $r\leftrightarrow\overline\one_c$. This is precisely the reason why we called the group~$\G$ a ``twin brother'' of~$\Hh$. In particular, by Proposition~\ref{prop:groupH} we conclude that $\overline\one_a\overline\one_c$ (corresponding to~$er\in\Hh$) has infinite order.

Similarly to the case of~$\Hh$, the group~$\hat\G/K_{\hat\G}$ acts faithfully on the tree~$T_{\mathcal O_{ac}}$, which is isomorphic to the tree~$Q_{\G}\cdot Y^*$ via isomorphisms $\psi_s$, $s\in Q_{\G}$. This isomorphism between trees induces a faithful action of~$\hat\G/K_{\hat\G}$ on~$Q_{\G}\cdot Y^*$ defined by the following wreath recursion:
\begin{equation}
\label{eqn_action_subtreeG}
\begin{array}{lcl}
\overline{\zero}&=&(\overline\one_a,\overline\zero_a,\overline\zero_c,\overline\zero_c)(a\,d\,b\,c),\\
\overline{\one}&=&(\overline\zero_a,\overline\one_a,\overline\one_c,\overline\one_c)(a\,d)(b\,c),
\end{array}
\end{equation}
where the sections $\overline\zero_a$, $\overline\zero_c$, $\overline\one_a$, and $\overline\one_c$ are defined by~\eqref{eqn_autom_sectionsG}. We will identify $\hat\G/K_{\hat\G}$ with the group defined by this wreath recursion. It follows from~\eqref{eqn_action_subtreeG} that
\[\overline{\zero}\cdot\overline{\one}=(\overline{\one}_a\overline{\one}_c, \overline{\zero}_a\overline{\one}_c, \overline{\zero}_c\overline{\zero}_a, \overline{\zero}_c\overline{\zero}_a)(c\,d).\]
Therefore, $\overline{\zero}\cdot\overline{\one}$ fixes vertex $a$ and its section $\overline{\one}_a\overline{\one}_c$ at this vertex has an infinite order by Proposition~\ref{prop:groupH}. Thus, $\overline{\zero}\cdot\overline{\one}$ itself has infinite order and $\hat\G/K_{\hat\G}$ is infinite.
\end{proof}

The proofs of the following corollary and lemma are identical to the proofs of Corollary~\ref{cor_stabs_diff} and Lemma~\ref{lem:transitivityH}.

\begin{corollary}
\label{cor_stabs_diffG}
The stabilizers of levels of~$T_{\mathcal O_{ac}}$ in~$\hat\G/K_{\hat\G}$ are pairwise different.
\end{corollary}

\begin{lemma}
\label{lem:transitivityG}
The action of~$\hat\G$ on the levels of~$T_{\mathcal O_{ac}}$ is transitive.
\end{lemma}

As in the case of $\Hh$, the last lemma can be rephrased in terms of the orbit tree of $\hat\G$ with the proof identical to the proof of Corollary~\ref{cor:path2}.

\begin{corollary}
\label{cor:path2G}
The orbit tree of~$\hat\G$ contains an infinite path from the root labeled by~$42^{\infty}$.
\end{corollary}

\begin{proof}[Proof of Proposition~\ref{prop:groupG}.]
The elements of~$\G$ corresponding to vertices $a(cb)^n$, $a(cb)^nc$, $n\geq0$ of odd and even levels of~$T_{\mathcal O_{ac}}$ respectively,
act nontrivially on the first level. By Lemma~\ref{lem:transitivityG} each word $v$ corresponding to a vertex in~$T_{\mathcal O_{ac}}$ can be moved by $\hat\G$ to one of these vertices. But this means that $v$ has a nontrivial section and thus must be a nontrivial element of~$\G$ itself.
\end{proof}

\subsection{Dependence of orbit automata on the representative of a symmetry class}
We conclude the paper with an example supporting the conclusion of Remark~\ref{rem:symmetric}. We will show that the relation between $\Hh$ and $\G$ shown in Figure~\ref{fig:G_Hh} is dependent on the symmetry class of the automata generating these groups. Moreover, we will construct a diagram showing this dependence in details.

\begin{example}
\label{ex:non_sym}
The automaton $\hat\C$ defined above in equalities~\eqref{eqn_autom_dual_defG} is symmetric to the automaton $\mathcal D$ (obtained from~$\hat\C$ by swapping letters $b$ and $d$ of the alphabet) with the following wreath recursion:
\begin{equation}
\label{eqn_autom_dual_defD}
\begin{array}{lcl}
\zero&=&(\one,\zero,\zero,\zero)(a\,b\,d\,c),\\
\one&=&(\zero,\one,\one,\one)(a\,b)(c\,d).\\
\end{array}
\end{equation}
However, the only nontrivial orbit automaton of the group~$\langle\hat\C\rangle$ is minimally symmetric to the automaton $\B$ defined by~\eqref{eqn_automH_def}, while the only nontrivial orbit automaton of~$\langle\D\rangle$ is minimally symmetric to the automaton $\C$ defined by~\eqref{eqn_autom_def} that is, in turn, not minimally symmetric to~$\B$.
\end{example}

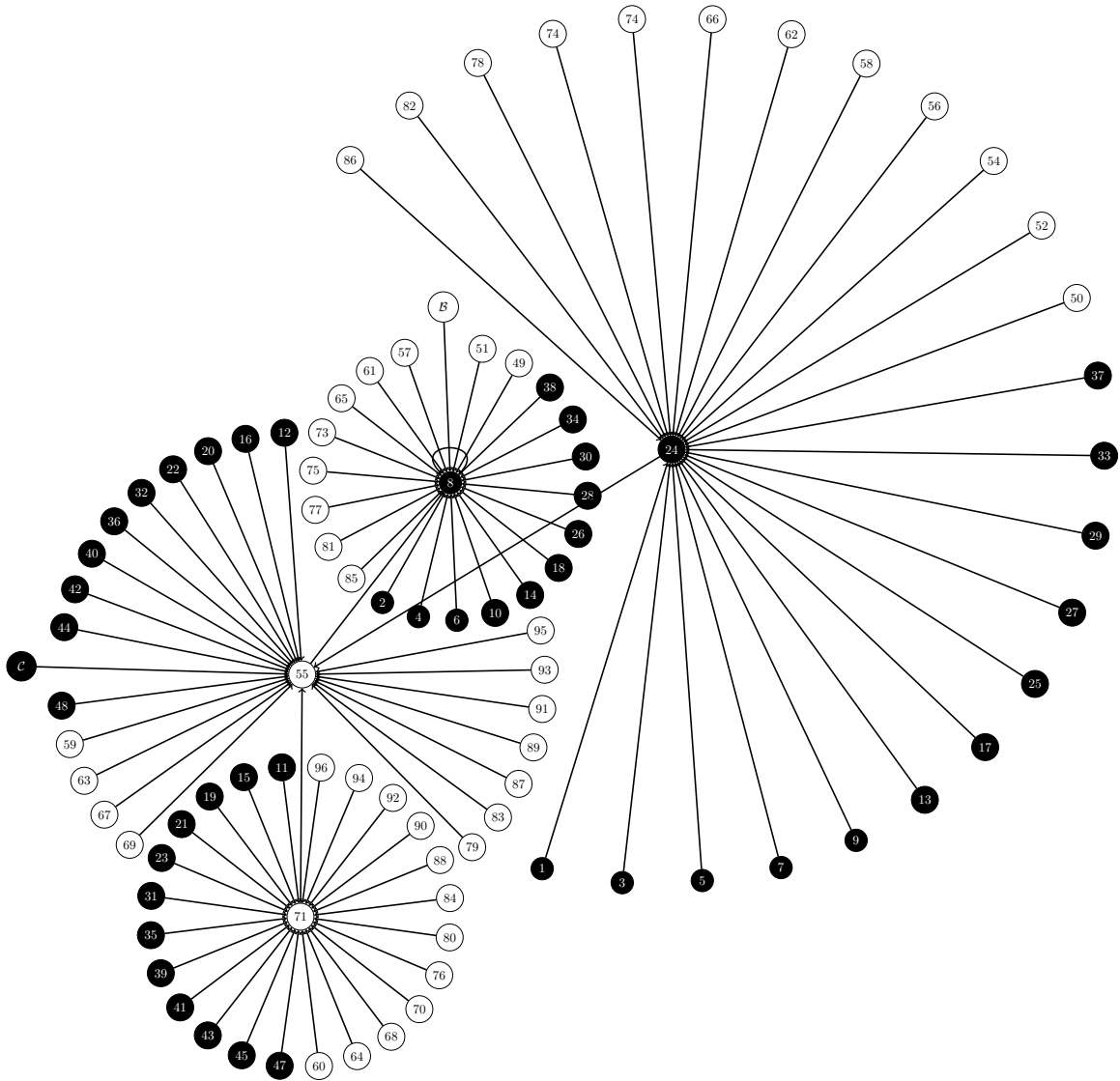
\begin{figure}
\begin{center}
\scalebox{.5}{\begin{tikzpicture}[scale=0.28,black]
  \begin{scope}[btn] 
  \node[btn] (86) at (2067.0bp,569.31bp) {1};
  \node[btn] (1) at (2290.0bp,530.31bp) {3};
  \node[btn] (2) at (1622.1bp,1311.7bp) {2};
  \node[btn] (3) at (2513.0bp,536.31bp) {5};
  \node[btn] (4) at (1722.2bp,1272.2bp) {4};
  \node[btn] (5) at (2732.4bp,572.5bp) {7};
  \node[btn] (6) at (1829.4bp,1262.1bp) {6};
  \node[btn] (7) at (2941.4bp,648.37bp) {9};
  \node[btn] (8) at (1811.4bp,1645.6bp) {8};
  \node[btn] (9) at (3132.9bp,761.35bp) {13};
  \node[btn] (10) at (1936.0bp,1282.5bp) {10};
  \node[btn] (11) at (1342.0bp,851.93bp) {11};
  \node[btn] (12) at (1349.0bp,1784.8bp) {12};
  \node[btn] (13) at (3301.0bp,908.26bp) {17};
  \node[btn] (14) at (2033.5bp,1332.5bp) {14};
  \node[btn] (15) at (1236.2bp,824.47bp) {15};
  \node[btn] (16) at (1240.6bp,1767.9bp) {16};
  \node[btn] (17) at (3439.7bp,1084.5bp) {25};
  \node[btn] (18) at (2112.9bp,1407.9bp) {18};
  \node[btn] (19) at (1140.9bp,770.45bp) {19};
  \node[btn] (20) at (1136.3bp,1733.6bp) {20};
  \node[btn] (21) at (1062.7bp,693.58bp) {21};
  \node[btn] (22) at (1038.9bp,1683.0bp) {22};
  \node[btn] (23) at (1007.2bp,599.12bp) {23};
  \node[btn] (24) at (2427.4bp,1738.4bp) {24};
  \node[btn] (25) at (3543.3bp,1283.3bp) {27};
  \node[btn] (26) at (2167.7bp,1502.7bp) {26};
  \node[btn] (27) at (3608.3bp,1497.9bp) {29};
  \node[btn] (28) at (2193.6bp,1609.2bp) {28};
  \node[btn] (29) at (3632.4bp,1720.9bp) {33};
  \node[btn] (30) at (2188.3bp,1718.6bp) {30};
  \node[btn] (31) at (978.0bp,493.49bp) {31};
  \node[btn] (32) at (950.99bp,1617.2bp) {32};
  \node[btn] (33) at (3614.8bp,1944.4bp) {37};
  \node[btn] (34) at (2152.4bp,1822.0bp) {34};
  \node[btn] (35) at (977.14bp,383.91bp) {35};
  \node[btn] (36) at (874.92bp,1538.0bp) {36};
  \node[btn] (38) at (2088.7bp,1911.1bp) {38};
  \node[btn] (39) at (1004.7bp,277.84bp) {39};
  \node[btn] (40) at (812.67bp,1447.6bp) {40};
  \node[btn] (41) at (1058.7bp,182.51bp) {41};
  \node[btn] (42) at (765.88bp,1348.3bp) {42};
  \node[btn] (43) at (1135.6bp,104.43bp) {43};
  \node[btn] (44) at (735.8bp,1242.7bp) {44};
  \node[btn] (45) at (1230.1bp,48.923bp) {45};
  \node[btn,inner sep=5pt,xshift=-30pt] (46) at (723.22bp,1133.7bp) {$\C$};
  \node[btn] (47) at (1335.8bp,19.766bp) {47};
  \node[btn] (48) at (728.46bp,1024.0bp) {48};
  \end{scope}
  \begin{scope}
  \node[ntb] (49) at (2002.5bp,1978.6bp)  {49};
  \node[ntb] (37) at (3556.1bp,2160.8bp) {50};
  \node[ntb] (50) at (3458.3bp,2362.6bp)  {52};
  \node[ntb] (51) at (1900.6bp,2019.0bp)  {51};
  \node[ntb] (52) at (3324.8bp,2542.8bp)  {54};
  \node[ntb,inner sep=5pt,yshift=30pt] (53) at (1791.6bp,2029.0bp)  {$\B$};
  \node[ntb] (54) at (3160.3bp,2695.1bp)  {56};
  \node[ntb] (55) at (1398.4bp,1111.0bp)   {55};
  \node[ntb] (56) at (2970.4bp,2814.3bp)  {58};
  \node[ntb] (57) at (1684.1bp,2007.8bp)  {57};
  \node[ntb] (58) at (2761.7bp,2896.3bp)  {62};
  \node[ntb] (59) at (751.4bp,916.65bp)    {59};
  \node[ntb] (60) at (1445.3bp,18.95bp)  {60};
  \node[ntb] (61) at (1587.0bp,1957.1bp)  {61};
  \node[ntb] (62) at (2541.4bp,2938.2bp)  {66};
  \node[ntb] (63) at (791.42bp,814.42bp)  {63};
  \node[ntb] (64) at (1551.4bp,46.531bp)  {64};
  \node[ntb] (65) at (1508.2bp,1881.1bp)  {65};
  \node[ntb] (66) at (2317.2bp,2938.5bp)  {74};
  \node[ntb] (67) at (847.47bp,720.03bp)  {67};
  \node[ntb] (68) at (1646.7bp,100.63bp)  {68};
  \node[ntb] (69) at (918.06bp,635.96bp)  {69};
  \node[ntb] (70) at (1724.8bp,177.55bp)  {70};
  \node[ntb] (71) at (1393.7bp,435.44bp)  {71};
  \node[ntb] (72) at (1780.2bp,272.06bp)  {76};
  \node[ntb] (73) at (1454.1bp,1785.9bp)  {73};
  \node[ntb] (74) at (2096.7bp,2897.3bp)  {74};
  \node[ntb] (75) at (1429.0bp,1679.3bp)  {75};
  \node[ntb] (76) at (1887.8bp,2816.0bp)  {78};
  \node[ntb] (77) at (1435.1bp,1569.9bp)  {77};
  \node[ntb] (78) at (1697.5bp,2697.4bp)  {82};
  \node[ntb] (79) at (1872.0bp,629.23bp)  {79};
  \node[ntb] (80) at (1809.3bp,377.71bp)  {80};
  \node[ntb] (81) at (1471.8bp,1466.7bp)  {81};
  \node[ntb] (82) at (1532.5bp,2545.6bp)  {86};
  \node[ntb] (83) at (1943.8bp,712.29bp)  {83};
  \node[ntb] (84) at (1810.1bp,487.3bp)    {84};
  \node[ntb] (85) at (1536.1bp,1378.1bp)  {85};
  \node[ntb] (87) at (2001.2bp,805.88bp)  {87};
  \node[ntb] (88) at (1782.5bp,593.34bp)  {88};
  \node[ntb] (89) at (2042.6bp,907.53bp)  {89};
  \node[ntb] (90) at (1728.4bp,688.63bp)  {90};
  \node[ntb] (91) at (2067.1bp,1014.6bp)  {91};
  \node[ntb] (92) at (1651.4bp,766.65bp)  {92};
  \node[ntb] (93) at (2073.9bp,1124.1bp)  {93};
  \node[ntb] (94) at (1556.9bp,822.08bp)  {94};
  \node[ntb] (95) at (2062.8bp,1233.3bp)  {95};
  \node[ntb] (96) at (1451.2bp,851.16bp)  {96};
  \end{scope}
  \begin{scope}[->,very thick]
  \path (8) edge[loop] (8)
  	(24) edge (55)
	(1) edge (24)
	(38) edge (8)
	(19) edge (71)
	(56) edge (24)
	(80) edge (71)
	(94) edge (71)
	(96) edge (71)
	(2) edge (8)
	(26) edge (8)
	(83) edge (55)
	(30) edge (8)
	(50) edge (24)
	(43) edge (71)
	(33) edge (24)
	(17) edge (24)
	(78) edge (24)
	(92) edge (71)
	(63) edge (55)
	(62) edge (24)
	(73) edge (8)
	(55) edge (8)
	(65) edge (8)
	(64) edge (71)
	(70) edge (71)
	(45) edge (71)
	(37) edge (24)
	(57) edge (8)
	(5) edge (24)
	(46) edge (55)
	(23) edge (71)
	(90) edge (71)
	(6) edge (8)
	(3) edge (24)
	(87) edge (55)
	(59) edge (55)
	(28) edge (8)
	(68) edge (71)
	(89) edge (55)
	(44) edge (55)
	(91) edge (55)
	(72) edge (71)
	(67) edge (55)
	(79) edge (55)
	(40) edge (55)
	(20) edge (55)
	(32) edge (55)
	(81) edge (8)
	(51) edge (8)
	(93) edge (55)
	(21) edge (71)
	(41) edge (71)
	(9) edge (24)
	(10) edge (8)
	(12) edge (55)
	(39) edge (71)
	(76) edge (24)
	(42) edge (55)
	(88) edge (71)
	(58) edge (24)
	(34) edge (8)
	(49) edge (8)
	(77) edge (8)
	(31) edge (71)
	(18) edge (8)
	(4) edge (8)
	(22) edge (55)
	(95) edge (55)
	(61) edge (8)
	(82) edge (24)
	(66) edge (24)
	(75) edge (8)
	(47) edge (71)
	(71) edge (55)
	(36) edge (55)
	(85) edge (8)
	(25) edge (24)
	(13) edge (24)
	(52) edge (24)
	(14) edge (8)
	(15) edge (71)
	(69) edge (55)
	(16) edge (55)
	(35) edge (71)
	(27) edge (24)
	(84) edge (71)
	(86) edge (24)
	(7) edge (24)
	(11) edge (71)
	(54) edge (24)
	(48) edge (55)
	(53) edge (8)
	(60) edge (71)
	(74) edge (24)
	(29) edge (24);
  \end{scope}
\end{tikzpicture}}
\end{center}
\caption{Relations between automata symmetric to~$\B$ (white circles) and those symmetric to~$\C$ (black circles).\label{fig-petals}}
\end{figure}

More precisely, there are 48 symmetric automata generating~$\G$ and 48 symmetric automata generating~$\Hh$. To each of these 96 automata one can apply the same operation of passing to the dual and then to (the only nontrivial) orbit automaton. In all cases one obtains an automaton symmetric to either~$\B$ or~$\C$. More precisely, half of 48 automata symmetric to~$\B$ produces automata symmetric to~$\C$ and the other half produces automata symmetric to~$\B$ again (so these automata are ``self-dual'' in this sense). And the same is true for automata symmetric to~$\C$. In Figure~\ref{fig-petals} we show the diagram describing precisely the above connections. The black vertices correspond to automata symmetric to~$\C$ and the white ones correspond to those symmetric to~$\B$. The arrows indicate that the target automaton was obtained from the source automaton by passing to the dual and then to the minimization of the nontrivial orbit automaton. Finally, the vertices are labeled by the numbers of corresponding automata in the lexicographic order on the set of all automata symmetric to~$\B$ and $\C$ respectively.

\section*{References}


\def\cprime{$'$} \def\cprime{$'$} \def\cprime{$'$}

\end{document}